\documentclass[11pt]{article}

\usepackage{times}

\usepackage{graphicx}

\usepackage{textcomp}
\usepackage{amsmath}
\usepackage{amsfonts}
\newcommand{\II}{\ensuremath{{\mathbb I}}}

\newtheorem{theorem}{Theorem}
\newtheorem{lemma}{Lemma}

\newtheorem{corollary}{Corollary}
\newtheorem{proposition}{Proposition}

\DeclareGraphicsExtensions{.eps, .jpg}

\title{ON BAYESIAN ``TESTIMATION'' AND ITS APPLICATION TO WAVELET THRESHOLDING}

\author{{\bf Felix Abramovich} \\
Department of Statistics and Operations Research, Tel Aviv
University
\and
{\bf Vadim Grinshtein} \\
Department of Mathematics, The Open University of Israel
\and
{\bf Athanasia Petsa \; and \; Theofanis Sapatinas} \\
Department of Mathematics and Statistics,
          University of Cyprus\\
}
\begin{document}
\maketitle

\date{}

%Figure inclusion macros

%Example: \fig{myfig}{An example of a good fit on worthless data}
%Inserts postscript figure myfig.ps into document, with caption given
%This is for .ps files not generated by transfig
%\fig{label=.ps prefix}{caption}

\begin{abstract}
We consider the problem of estimating the unknown response function
in the Gaussian white noise model. We first utilize the
recently developed Bayesian maximum a posteriori ``testimation''
procedure of Abramovich et al. (2007) for recovering an unknown high-dimensional Gaussian mean
vector. The existing results for its upper error bounds over various
sparse $l_p$-balls are extended to more general cases. We show that,
for a properly chosen prior on the number of non-zero entries of the
mean vector, the corresponding adaptive estimator is asymptotically
minimax in a wide range of sparse and dense $l_p$-balls.

The proposed procedure is then applied in a wavelet context
to derive adaptive global and level-wise wavelet estimators of
the unknown response function in the Gaussian white noise
model. These estimators are then proven to be, respectively,
asymptotically near-minimax and minimax in a wide range of Besov balls. 
These results are also extended to the estimation of derivatives of the response
function.

Simulated examples are conducted to illustrate the performance of
the proposed level-wise wavelet estimator in finite
sample situations, and to compare it with several existing
counterparts.
\end{abstract}

%\begin{keywords}
%Adaptive estimation; Besov
%Spaces; Gaussian sequence model; Gaussian white noise model; $l_{p}$-balls;
%Multiple testing; Thresholding; Wavelet estimation.
%\end{keywords}

\newpage
\section{Introduction}\label{In2}
We consider the problem of estimating the unknown response function
in the Gaussian white noise model, where
one observes Gaussian processes $Y_n(t)$ governed by
\begin{equation}
\label{eq:FaGauWh} dY_n(t) = f(t)dt + \frac{\sigma}{\surd{n}}\, dW(t),
\quad t \in [0,1].
\end{equation}
The noise parameter $\sigma >0$ is assumed to be known, $W$ is a
standard Wiener process, and $f \in L^2[0,1]$ is the unknown
response function. Under some smoothness constraints on $f$, such a
model is asymptotically equivalent in Le Cam sense to the standard
nonparametric regression setting (Brown \& Low, 1996).

In a consistent estimation theory, it is well-known that
$f$ should possess some smoothness properties. We assume that $f$
belongs to a Besov ball $B_{p,q}^{s}(M)$ of a radius $M>0$, where
$0 < p,q \leq \infty$ and $s > \max(0,1/p-1/2)$. The latter restriction
ensures that the corresponding Besov spaces are embedded in
$L^2[0,1]$. The parameter $s$ measures the degree of smoothness while
$p$ and $q$ specify the type of norm used to measure the smoothness.
Besov classes contain various traditional smoothness spaces such as
H\"older and Sobolev spaces as special cases. However, they also
include different types of spatially inhomogeneous functions (Meyer, 1992).

The fact that wavelet series constitute unconditional bases for Besov spaces
has caused various wavelet-based estimation procedures to be widely
used for estimating the unknown response $f \in B^{s}_{p,q}(M)$ in
the Gaussian white noise model (\ref{eq:FaGauWh}). 
The standard wavelet approach for the
estimation of $f$ is based on finding the empirical wavelet
coefficients of the data and denoising them, usually by
some type of a thresholding rule. Transforming them back to the function
space then yields the resulting estimate. The main statistical
challenge in such an approach is a proper choice of a thresholding
rule. A series of various wavelet thresholds originated by different
ideas has been proposed in the literature during the last decade,
e.g., the universal threshold (Donoho \& Johnstone, 1994a), Stein's
unbiased risk estimation
threshold (Donoho \& Johnstone, 1995), 
the false discovery rate threshold (Abramovich \&
Benjamini, 1996), the cross-validation threshold (Nason, 1996), the Bayes
threshold (Abramovich et al., 1998) and the empirical
Bayes threshold (Johnstone \& Silverman, 2005).

Abramovich \& Benjamini (1996) demonstrated that thresholding can be
viewed as a multiple hypothesis testing procedure, where one first
simultaneously tests the wavelet coefficients of the unknown
response function, for significance. The coefficients concluded to be
significant are then estimated by the corresponding empirical
wavelet coefficients of the data, while the non-significant ones are discarded. 
Such a ``testimation'' procedure evidently mimics a hard thresholding
rule. Various choices for adjustment to multiplicity on the testing
step lead to different thresholds. In particular, the universal
threshold of Donoho \& Johnstone (1994a) and the false discovery rate threshold 
of Abramovich \& Benjamini (1996) fall within such a framework
corresponding to Bonferroni and false discovery rate multiplicity corrections,
respectively.

In this paper, we proceed along the lines of ``testimation''
approach, where we utilize the recently developed maximum a
posteriori Bayesian multiple testing procedure of Abramovich
\& Angelini (2006). Their hierarchical prior model
is based on imposing a prior distribution on the number of
false null hypotheses. Abramovich et al. (2007)
applied this approach to estimating a high-dimensional Gaussian mean vector and
showed its minimax optimality 
where the unknown mean vector was assumed to be sparse.

We first extend the results of Abramovich et al. (2007) to more general settings.
Consider the problem of estimating
an unknown high-dimensional Gaussian mean vector, where one observes
$y_i$ governed by
\begin{equation}
y_i=\mu_i+\sigma_n\, z_i,\quad i=1,2,\ldots,n. \label{eq:sequence}
\end{equation}
The variance $\sigma^2_n >0$, that may depend on $n$, is assumed to
be known, $z_i$ are independent $N(0,1)$ random variables, and the
unknown mean vector $\mu=(\mu_1,\ldots,\mu_n)'$ is assumed to lie in
a strong $l_p$-ball $l_p[\eta_n]$, $0 < p \leq \infty$, of a normalized
radius $\eta_n$, that is, $||\mu||_p \leq C_n$, where
$C_n=n^{1/p}\sigma_n \eta_n$.
Abramovich et al. (2007)
considered  the Gaussian sequence model (\ref{eq:sequence})
with $\sigma_n^2=\sigma^2$ and derived upper error bounds for the
quadratic risk of an adaptive Bayesian maximum a posteriori estimator of $\mu$ 
in the sparse
case, where $0< p<2$ and $\eta_n \rightarrow 0$ as $n \rightarrow
\infty$. We extend their results for all combinations of $p$ and
$\eta_n$ and for the variance in (\ref{eq:sequence}) that may depend 
on $n$. We show, in particular, that for a properly chosen prior
distribution on the number of non-zero entries of $\mu$, the
corresponding estimator, up to a constant factor,
is asymptotically minimax for almost all $l_p$-balls
including both sparse and dense cases.

We then apply the proposed approach to the wavelet
thresholding estimation in the Gaussian white noise model
(\ref{eq:FaGauWh}). We show that, under mild conditions on the prior
distribution on the number of non-zero wavelet coefficients, the resulting
global wavelet estimator of $f$, up to a logarithmic factor, attains the minimax
convergence rates simultaneously over the entire range of Besov balls.
Furthermore, we demonstrate that estimating wavelet coefficients
at each resolution level separately, allows one to remove
the extra logarithmic factor. Moreover, the procedure can also be
extended to the estimation of derivatives of $f$. These results, 
in some sense, complement the adaptively minimax empirical Bayes estimators of
Johnstone \& Silverman (2005).

\section{Estimation in the Gaussian sequence model} \label{sec:MAP}
\subsection{Bayesian maximum a posteriori estimation procedure} \label{subsec:test}
We start with reviewing the Bayesian maximum a posteriori estimation procedure 
for the Gaussian sequence model (\ref{eq:sequence}) developed by
Abramovich et al. (2007).

For this model, consider the multiple hypothesis testing problem,
where we wish to simultaneously test
$$H_{0i}:\mu_i=0 \quad \text{versus} \quad H_{1i}: \mu_i \ne 0, \quad i=1,2,\ldots,n.
$$
A configuration of true and false null hypotheses is uniquely
defined by the indicator vector $x=(x_1,...,x_n)'$, where $x_i=\II(\mu_i \ne
0)$ and $\II(A)$ denotes the indicator
function of the set $A$. Let $\kappa=x_1+...+x_n=||\mu||_0$ be the
number of non-zero $\mu_i$, i.e., $||\mu||_0
=\#\{i:\,\mu_i \ne 0\}$. Assume some prior distribution $\pi_n$ on
$\kappa$ with $\pi_n(\kappa)>0,\; \kappa=0,\ldots,n$. For a given
$\kappa$, all the corresponding different vectors $x$ are assumed 
to be equally likely a priori, that is,
conditionally on $\kappa$,
$$
{\rm pr} \bigg(x \mid \sum_{i=1}^n x_i=\kappa \bigg)={n \choose
\kappa}^{-1}.
%pr\bigg(x \,\big|\, \sum_{i=1}^n x_i=\kappa \bigg)={n \choose
%\kappa}^{-1}.
$$
Naturally, $\mu_i \mid x_i=0 \sim \delta_0$, where $\delta_0$ is
a probability atom at zero. To complete the prior specification, we
assume that $\mu_i \mid x_i=1 \sim N(0,\tau_n^2)$.

For the proposed hierarchical prior, the posterior probability of a
given vector $x$ with $\kappa$ non-zero entries is
\begin{equation}
\pi_n(x,\kappa \mid y) \propto {n \choose
\kappa}^{-1}\pi_n(\kappa)\,\II\bigg(\sum_{i=1}^n x_i=\kappa \bigg)
\prod_{i=1}^n (B_i^{-1})^{x_i}, \label{eq:post}
\end{equation}
where the Bayes factor $B_i$ of $H_{0i}$ is
\begin{equation}
B_i=\surd(1+\gamma_n)\; \exp\left\{-\frac{y_i^2}{2\sigma_n^2(1+1/\gamma_n)}\right\}
\label{eq:bf}
\end{equation}
and $\gamma_n=\tau_n^2/\sigma_n^2$ is the variance ratio
(Abramovich \& Angelini, 2006).

Given the posterior distribution $\pi_n(x,\kappa \mid y)$, we apply
the maximum a posteriori rule to choose the most likely indicator vector. Generally, to
find the posterior mode of $\pi_n(x,\kappa \mid y)$, one should look
through all $2^n$ possible sequences of zeroes and ones. However,
for the proposed model, the number of candidates for a mode is, in
fact, reduced to $n+1$ only. Indeed, let $\hat{x}(\kappa)$ be a
maximizer of (\ref{eq:post}) for a fixed $\kappa$ that indicates the
most plausible vector $x$ with $\kappa$ non-zero entries. From
(\ref{eq:post}), it follows immediately that $\hat{x}_i(\kappa)=1$
at the $\kappa$ entries corresponding to the smallest Bayes factors
$B_i$ and zeroes otherwise. Due to the monotonicity of $B_i$ in
$|y|_i$ in (\ref{eq:bf}), it is equivalent to setting
$\hat{x}_i(\kappa)=1$ for the $\kappa$ largest $|y|_i$
and zeroes for others. The proposed Bayesian multiple testing procedure
then leads to finding $\hat{\kappa}$ that maximizes
$$
\log \pi_n(\hat{x}(\kappa),\kappa \mid y) = c + \sum_{i=1}^{\kappa}
y_{(i)}^2+2\sigma_n^2(1+1/\gamma_n)\log\left\{{n \choose
\kappa}^{-1}\pi_n(\kappa) (1+\gamma_n)^{-\kappa/2}\right\}
$$
for some constant $c$ 
or, equivalently, minimizes
$$
\sum_{i=\kappa+1}^n y_{(i)}^2+2\sigma_n^2(1+1/\gamma_n)\log\left\{{n
\choose \kappa}\pi^{-1}_n(\kappa)
(1+\gamma_n)^{\kappa/2}\right\},
$$
where $|y|_{(1)} \geq \ldots \geq |y|_{(n)}$. The ${\hat \kappa}$
null hypotheses corresponding to
$|y|_{(1)},\ldots,|y|_{(\hat{\kappa})}$ are rejected. The resulting
Bayesian estimation yields a hard thresholding with a threshold
$\hat{\lambda}_{\rm MAP}=|y|_{(\hat{\kappa})}$, i.e.,
\begin{equation}
\hat{\mu}_i=\left\{
\begin{array}{ll}
y_i, & |y_i| \geq \hat{\lambda}_{\rm MAP}, \\
0, & {\rm otherwise.}
\end{array}
\right.
\label{eq:MAP}
\end{equation}
If $\hat{\kappa}=0$, then all $y_i$, $i=1,2,\ldots,n$, are
thresholded and $\hat{\mu} \equiv 0$.

From a frequentist view, the above estimator 
$\hat{\mu}=(\hat{\mu}_1,\ldots,\hat{\mu}_n)'$ in  (\ref{eq:MAP}) is evidently a penalized
likelihood estimator with the complexity penalty
\begin{equation}
P_n(\kappa)=2\sigma_n^2(1+1/\gamma_n)\log\left\{{n \choose
\kappa}\pi^{-1}_n(\kappa) (1+\gamma_n)^{\kappa/2}\right\}.
\label{eq:penalty}
\end{equation}
In this sense, it can be also considered within the
framework of Birg\'e \& Massart (2001).
We will discuss these relations in the following section in more details.

\subsection{Upper error bounds}
\label{subsec:bounds} Abramovich et al. (2007, Theorem 6)
obtained upper error bounds for the $l^2$-risk of (\ref{eq:MAP})
in the Gaussian sequence model
(\ref{eq:sequence}) for sparse $l_p[\eta_n]$-balls, where
$0<p<2$ and $\eta_n \rightarrow 0$ as $n \rightarrow \infty$. We
extend now these results to more general settings.

Fix a prior distribution $\pi_n(\kappa)>0,\;\kappa=0,\ldots,n$, on
the number of non-zero entries of $\mu$, and let
$\gamma_n=\tau_n^2/\sigma_n^2$ be the variance ratio.

\begin{proposition} \label{prop:prop1}
Let $\hat{\mu}$ be the estimator (\ref{eq:MAP}) of $\mu$ in
the Gaussian sequence model (\ref{eq:sequence}), where $\mu \in
l_p[\eta_n]$, $0 < p \leq \infty$.
Assume that there exist positive
constants $\gamma_{-}$ and $\gamma_{+}$ such that
$\gamma_{-} \leq \gamma_n \leq \gamma_{+}$.
\begin{enumerate}
\item Let $0 < p \leq \infty$. 
Assume that $\pi_n(n) \geq e^{-c_0 n}$ for some $c_0>0$.
Then, as $n \rightarrow \infty$,
$$
\sup_{\mu \in l_p[\eta_n]} E(||\hat{\mu}-\mu||^2_2)=O(n\sigma_n^2).
$$

\item Let $2 \leq p \leq \infty$. Assume that there exists $\beta \geq 0$ such that
$\pi_n(0) \geq n^{-c_1 n^{-\beta}}$ for some $c_1>0$. Then, as $n \rightarrow \infty$,
$$
\sup_{\mu \in l_p[\eta_n]} E(||\hat{\mu}-\mu||^2_2) =
O(\sigma_n^2n\eta_n^2)+O(\sigma_n^2n^{-\beta}\log n). 
$$

\item Let $0 < p < 2$. Assume $\pi_n(\kappa) \geq (\kappa/n)^{c_2 \kappa}$ for all $\kappa=1,2,\ldots,\alpha_n n$,
where $n^{-1}(2 \log n)^{p/2} \leq \alpha_n \leq \exp\{-c(\gamma_n)\},\;c(\gamma_n)=8(\gamma_n+3/4)^2 > 9/2$,
and for some $c_2>0$. Then, as $n \rightarrow \infty$,
$$
\sup_{\mu \in l_p[\eta_n]} E(||\hat{\mu}-\mu||^2_2)=
O\big\{\sigma_n^2 n \eta_n^p (2\log \eta_n^{-p})^{1-p/2}\big\}
$$
for all $n^{-1}(2 \log n)^{p/2} \leq \eta^p_n \leq \alpha_n$.

\item Let $0 < p < 2$. Assume that there exists $\beta \geq 0$ such that
$\pi_n(0) \geq n^{-c_1n^{-\beta}}$ for some $c_1>0$. Then, as $n \rightarrow \infty$,
$$
\sup_{\mu \in l_p[\eta_n]} E(||\hat{\mu}-\mu||^2_2) = O(\sigma_n^2
n^{2/p} \eta^2_n) +O(\sigma_n^2 n^{-\beta} \log n) 
$$
for all $\eta^p_n < n^{-1}(2 \log n)^{p/2}$.
\end{enumerate}
\end{proposition}

The proof of Proposition \ref{prop:prop1} is given in the Appendix.
Similar to Abramovich et al. (2007),
analogous results can be obtained for other types of balls, e.g.,
weak $l_p$-balls, $0 < p < \infty$, and $l_0$-balls, with necessary changes in the
proofs (Petsa, 2009, Chapter 3).

Since the prior assumptions in Proposition \ref{prop:prop1} do not
depend on the parameters $p$ and $\eta_n$ of the $l_p$-ball, the
estimator (\ref{eq:MAP}) is inherently adaptive. The condition on
$\pi_n(n)$ guarantees that its risk is always
bounded by an order of $n\sigma^2_n$, corresponding to the risk of
the maximum likelihood estimator, $\hat{\mu}^{MLE}_i=y_i$, in the
Gaussian sequence model (\ref{eq:sequence}).

The following corollary of Proposition \ref{prop:prop1} essentially
defines dense and sparse zones for $2 \leq p \leq \infty$, and 
dense, sparse and super-sparse zones for $0 < p < 2$ of different
behavior for the quadratic risk of the proposed estimator (\ref{eq:MAP}). 
To evaluate its accuracy, we also compare the resulting
risks with the corresponding minimax risks
$R(l_p[\eta_n])=\inf_{\tilde{\mu}} \sup_{\mu \in l_p[\eta_n]}
E(||\tilde{\mu}-\mu||^2_2)$ that can be found, e.g., in Donoho \&
Johnstone (1994b). In what follows 
$g_1(n) \asymp g_2(n)$ denotes $0 < \liminf\{g_1(n)/g_2(n)\}
\leq \limsup\{g_1(n)/g_2(n)\} < \infty$ as $n \rightarrow \infty$.

\begin{corollary} \label{col:col1}
Let $\hat{\mu}$ be the estimator (\ref{eq:MAP}) of $\mu$ in
the Gaussian sequence model (\ref{eq:sequence}), where $\mu \in
l_p[\eta_n]$, $0 < p \leq \infty$.
Assume that there exist positive
constants $\gamma_{-}$ and $\gamma_{+}$ such that
$\gamma_{-} \leq \gamma_n \leq \gamma_{+}$. Define
$c(\gamma_n)=8(\gamma_n+3/4)^2  > 9/2$
and let the prior $\pi_n$ satisfy the following conditions:
\begin{enumerate}
\item $\pi_n(0) \geq n^{-c_1 n^{-\beta}}$ for some $\beta \geq 0$ and
$c_1>0$;
\item $\pi_n(\kappa) \geq (\kappa/n)^{c_2 \kappa}$ for all
$\kappa=1,2,\ldots,\alpha n$, where $\alpha=\exp(-9/2)$ or
$\alpha=\exp\{-c(\gamma_{-})\}$ if $\gamma_{-}$ is
known, and for some $c_2>0$;
\item $\pi_n(n) \geq e^{-c_0 n}$ for some $c_0>0$. 
\end{enumerate}
Then, as $n \rightarrow \infty$, depending on $p$ and $\eta_n$, one
has:
\begin{description}
\item[]Case 1. Let $\;0 < p \leq \infty,\;\eta_n^p > \alpha$. Then,
$$
\sup_{\mu \in l_p[\eta_n]} E(||\hat{\mu}-\mu||^2_2)=O(n\sigma_n^2),
\;\;\; R(l_p[\eta_n]) \asymp n \sigma_n^2
$$

\item[]Case 2. Let $\;2 \leq p \leq \infty,\;\eta_n^p \leq \alpha$. Then,
$$
\sup_{\mu \in l_p[\eta_n]} E(||\hat{\mu}-\mu||^2_2)=
O(\sigma_n^2n\eta_n^2)+O(\sigma_n^2n^{-\beta}\log n),\;\;\;
R(l_p[\eta_n]) \asymp \sigma_n^2 n \eta_n^2
$$

\item[]Case 3. Let $\;0<p<2,\;n^{-1}(2 \log n)^{p/2} \leq \eta^p_n \leq
\alpha$. Then,
$$
\sup_{\mu \in l_p[\eta_n]} E(||\hat{\mu}-\mu||^2_2)=
O\big\{\sigma_n^2 n \eta_n^p (2\log \eta_n^{-p})^{1-p/2}\big\},\;\;\;
R(l_p[\eta_n]) \asymp \sigma_n^2 n \eta_n^p(2\log \eta_n^{-p})^{1-p/2}
$$

\item[]Case 4. Let $\;0<p<2,\;\eta^p_n < n^{-1}(2 \log n)^{p/2}$. Then,
$$
\sup_{\mu \in l_p[\eta_n]} E(||\hat{\mu}-\mu||^2_2) = O(\sigma_n^2
n^{2/p} \eta^2_n) +O(\sigma_n^2 n^{-\beta} \log n),\;\;\;
R(l_p[\eta_n]) \asymp \sigma_n^2 n^{2/p}\eta_n^2
$$

\end{description}

\end{corollary}
For $\beta=0$ one can easily verify that all three conditions
of Corollary \ref{col:col1} are satisfied, for example, for the
truncated geometric prior ${\rm TrGeom}(1-q),\;0 < q < 1$, where
$\pi_n(\kappa)=(1-q)q^\kappa/(1-q^{n+1}),\;\kappa=0,\ldots,n$.
On the other hand, for any $\beta$, no binomial
prior ${\rm Bin}(n,p_n)$ can ``kill three birds with one stone''.
The requirement $\pi_n(0)=(1-p_n)^n \geq n^{-c_1n^{-\beta}}$ necessarily
implies $p_n \rightarrow 0$ as $n \rightarrow \infty$. However, to
satisfy $\pi_n(n)=p_n^n \geq e^{-c_0 n}$, one needs $p_n \geq e^{-c_0}$. 

The impact of Corollary \ref{col:col1} is that, up to a constant multiplier,
the proposed estimator (\ref{eq:MAP}) is adaptively minimax for almost all
$l_p$-balls, $0 < p \leq \infty$, except those with very small
normalized radiuses, where
$\eta_n^2=o(n^{-(\beta+2/\min(p,2))}\log n)$. Hence, while the
optimality of most the existing threshold estimators, e.g.,
universal, Stein's unbiased risk, false discovey rate, has been established only
over various sparse settings, the Bayesian estimator (\ref{eq:MAP}) is 
appropriate for both sparse and
dense cases. To the best of our knowledge, such a wide adaptivity
range can be compared only with the penalized likelihood estimators
of Birg\'e \& Massart (2001) and the empirical Bayes threshold
estimators of Johnstone \& Silverman (2004b, 2005); see Petsa (2009, Chapter 3) 
for more details.

In fact, as we have mentioned, there are interesting asymptotic relationships 
between the Bayesian estimator (\ref{eq:MAP}) and the penalized likelihood
estimator of Birg\'e \& Massart (2001) that may explain their
similar behavior. For estimating the normal mean vector in
(\ref{eq:sequence}) within $l_p$-balls, Birg\'e \& Massart (2001)
considered a penalized likelihood estimator with a specific
complexity penalty
\begin{equation}
\label{eq:penBM-1}\tilde{P}_n(\kappa)=C \sigma_n^2 \kappa \{1+\surd{(2
L_\kappa)}\}^2,
\end{equation}
where $L_\kappa=\log(n/\kappa)+(1+\theta)(1+\log(n)/\kappa)$ for
fixed $C>1$ and $\theta > 0$ (Birg\'e \& Massart, 2001, Section 6.3). For
large $n$ and $\kappa< n/e$, this penalty is approximately of the
following form:
\begin{equation}
\label{eq:penBM-2} \tilde{P}_n(\kappa) \sim 2 \sigma_n^2 c \kappa L_\kappa
\sim 2 \sigma_n^2 \tilde{c}_1 \left\{\log{n \choose \kappa} +
\tilde{c}_2 \kappa\right\}
\end{equation}
for some positive constants $c,\;\tilde{c}_1,\;\tilde{c}_2 > 1$; see
also Lemma \ref{lem:lem1} in the Appendix. Thus, within this range,
$\tilde{P}_n$ in (\ref{eq:penBM-1})-(\ref{eq:penBM-2}) behaves in a way
similar to a particular case of the penalty $P_n$ in
(\ref{eq:penalty}) corresponding to the geometric type prior
$\pi_n(\kappa) \propto (1/\tilde{c}_2)^\kappa$. This prior
satisfies the second condition on $\pi_n$ of Corollary
\ref{col:col1}. Such a Bayesian interpretation can also be helpful
in providing some intuition behind the penalty $\tilde{P}_n$
motivated in Birg\'e \& Massart (2001) mostly due technical reasons.
In addition, under the conditions of Corollary \ref{col:col1},
$P_n(n) \sim \tilde{P}_n(n) \sim c n$.

Furthermore, for sparse cases, where $\kappa \ll n$, under the conditions
on the prior $\pi_n$ of Corollary \ref{col:col1}, both penalties
$P_n$ and $\tilde{P}_n$ 
are of the same so-called $2\kappa\log(n/\kappa)$-type
penalties of the form $2\sigma_n^2\zeta \kappa\{\log(n/\kappa)+
c_{\kappa,n}\},$ where $\zeta>1$ and $c_{\kappa,n}$ is negligible
relative to $\log(n/\kappa)$. Such type of penalties has appeared within
different frameworks in a series of recent works on estimation and
model selection (Foster \& Stine, 1999; George \& Foster,
2000; Birg\'e \& Massart, 2001; Abramovich et al., 2006;
Abramovich et al., 2007).

\medskip

\section{Bayesian maximum a posteriori wavelet estimation in the Gaussian white noise model} \label{sec:MAPwavelet}
\subsection{General algorithm} \label{subsec:general}
In this section we apply the results of Section \ref{sec:MAP} on
estimation in the Gaussian sequence model (\ref{eq:sequence})
to wavelet estimation of the unknown response function $f$ in the
Gaussian white noise model (\ref{eq:FaGauWh}).

Given a compactly supported scaling function $\phi$ of regularity
$r>s$ and the corresponding mother wavelet $\psi$, one can generate
an orthonormal wavelet basis on the unit interval from a finite
number $C_{j_0}$ of scaling functions $\phi_{j_0k}$ at a primary
resolution level $j_0$ and wavelets $\psi_{jk}$ at resolution levels
$j \geq j_0$ and scales $k=0,\ldots,2^j-1$ (Cohen et al., 1993;
Johnstone \& Silverman, 2004a). For
clarity of exposition, we use the same notation for interior and
edge wavelets, and in what follows denote $\phi_{j_0k}$ by
$\psi_{j_0-1,k}$.

Then, $f$ is expanded in the orthonormal wavelet series on
$[0,1]$ as
$$
f(t)=\sum_{j=j_0-1}^\infty \sum_{k=0}^{2^j-1} \theta_{jk} \psi_{jk}(t),
$$
where $\theta_{jk}=\int_0^1 f(t)\psi_{jk}(t)dt$. In the wavelet
domain, the Gaussian white noise model (\ref{eq:FaGauWh})
becomes
$$
Y_{jk}=\theta_{jk}+\epsilon_{jk},\;\;\;j \geq
j_0-1,\;k=0,\ldots,2^j-1,
$$
where the empirical wavelet coefficients $Y_{jk}$ are given by
$Y_{jk}=\int_0^1 \psi_{jk}(t)dY(t)$ and $\epsilon_{jk}$ are
independent $N(0,\sigma^2/n)$ random variables.

Define $J=\log_2 n$. Estimate wavelet coefficients $\theta_{jk}$ at
different resolution levels $j$ by the following scheme:
\begin{enumerate}
\item set $\hat{\theta}_{j_0-1,k}=Y_{j_0-1,k}$;
\item apply the Bayesian estimation procedure of Abramovich et al. (2007)
described in Section \ref{sec:MAP} to estimate
$\theta_{jk}$ at resolution levels $j_0 \leq j < J$ by the
corresponding $\hat{\theta}_{j,k}$;
\item set $\hat{\theta}_{jk}=0,\;j \geq J$.
\end{enumerate}
The resulting wavelet estimator $\hat{f}_n$ of $f$ is then
defined as
\begin{equation}
\hat{f}_n(t)=\sum_{k=0}^{C_{j_0}-1}Y_{j_0-1,k}\psi_{j_0-1,k}(t)+
\sum_{j=j_0}^{J-1} \sum_{k=0}^{2^j-1} \hat{\theta}_{jk}
\psi_{jk}(t).
\label{eq:waveMAP}
\end{equation}

Theorem \ref{th:global} below shows that, under mild conditions on the prior
$\pi_n$, the resulting global wavelet estimator (\ref{eq:waveMAP}) of $f$, 
where the estimation procedure is applied to the entire set of
wavelet coefficients at all resolution levels $j_0 \leq j < J$, up
to a logarithmic factor, attains the minimax convergence rates
over the whole range of Besov classes. 
Furthermore, Theorem \ref{th:level0} demonstrates that performing the 
estimation procedure at each resolution level separately allows one to remove
the extra logarithmic factor.
Moreover, a level-wise version of (\ref{eq:waveMAP}) 
allows one to estimate the derivatives of $f$ at optimal convergence
rates as well.

\subsection{Global wavelet estimator} \label{subsec:global}
The number of wavelet coefficients at all resolution levels up to
$J$ is $\tilde{n}=2^J-2^{j_0} \sim n$ for large $n$. Let
$\pi_n(\kappa)>0,\;\kappa=0,\ldots,\tilde{n}$, be a prior distribution
on the number of non-zero wavelet coefficients of $f$ at all resolution
levels $j_0 \leq j < J$, and let the prior variance of non-zero
coefficients at the $j$th resolution level be $\tau_j^2/n$; the
corresponding level-wise variance ratios are
$\gamma_j=\tau_j^2/\sigma^2$.

It is well-known (Donoho \& Johnstone, 1998) that, as $n \rightarrow \infty$,
the minimax convergence rate for the $L^2$-risk of estimating the unknown
response function $f$ in the model
(\ref{eq:FaGauWh}) over Besov balls $B^{s}_{p,q}(M)$, where
$0 < p,q \leq \infty,\;s>\max(0,1/p-1/2)$ and $M>0$, is given by
$$
\inf_{\tilde{f}_n}\sup_{f \in B^s_{p,q}(M)} E(||\tilde{f}_n-f||_2^2)
\asymp n^{-2s/(2s+1)}.
$$
%where $g_1(n) \asymp g_2(n)$ denotes $0 < \liminf\{g_1(n)/g_2(n)\}
%\leq \limsup\{g_1(n)/g_2(n)\} < \infty$ as $n \rightarrow \infty$.

\medskip

\begin{theorem} \label{th:global}
Let $\psi$ be a mother wavelet of regularity $r$ and let $\hat{f}_n$
be the corresponding global wavelet estimator (\ref{eq:waveMAP}) of $f$ in the
Gaussian white noise model (\ref{eq:FaGauWh}), where $f \in
B^s_{p,q}(M)$, $0 < p,q \leq \infty$, $1/p < s < r$
and $M>0$. Assume that there exist positive constants
$\gamma_{-}$ and $\gamma_{+}$ such that
$\gamma_{-} \leq \gamma_j \leq \gamma_{+}$ for all
$j=j_0,\ldots,J-1$. Let the prior $\pi_n$ satisfy
$\pi_n(\kappa) \geq (\kappa/n)^{c \kappa}$ for all
$\kappa=1,2,\ldots,\exp(-9/2)n$ or, for a shorter range
$\kappa=1,2,\ldots,\exp\{-c(\gamma_{-})\}n$ if
$\gamma_{-}$ is known. Then, as $n \rightarrow \infty$,
\begin{equation}
\sup_{f \in B^s_{p,q}(M)} E(||\hat{f}_n-f||_2^2)=
O\Bigg\{\bigg(\frac{\log n}{n}\bigg)^{\frac{2s}{2s+1}}\Bigg\} .
\label{eq:global}
\end{equation}
\end{theorem}

\medskip

The proof of Theorem \ref{th:global} is based on the
relationship between the smoothness conditions on functions within Besov
spaces and the conditions on their wavelet coefficients. Namely, if
$f \in B^s_{p,q}(M)$, then the sequence of its wavelet coefficients
$\{\theta_{jk},\;k=0,\ldots,2^j-1,\;j=j_0,\ldots,J-1\}$
belongs to a weak $l_{2/(2s+1)}$-ball of a radius $a M$, where the
constant $a$ depends only on a chosen wavelet basis (Donoho, 1993, Lemma 2).
One can then apply the corresponding results of
Abramovich et al. (2007) for estimation over
weak $l_p$-balls. Details of the proof of Theorem \ref{th:global}
are given in the Appendix.

The resulting global wavelet estimator does not rely on the
knowledge of the parameters $s$, $p$, $q$ and $M$ of a specific
Besov ball and it is, therefore, inherently adaptive. Theorem
\ref{th:global} establishes the upper bound for its $L^2$-risk and
shows that the resulting adaptive global wavelet estimator is
asymptotically near-optimal within the entire range of Besov balls.
In fact, the additional logarithmic factor in (\ref{eq:global}) is
the unavoidable minimal price for adaptivity for any global wavelet
threshold estimator (Donoho et al., 1995; Cai,
1999), and in this sense, the upper bound for the convergence rates
in (\ref{eq:global}) is sharp. To remove this logarithmic factor one
should consider level-wise thresholding. 

\subsection{Level-wise wavelet estimator} \label{subsec:level}
Consider now the level-wise version of the wavelet estimator (\ref{eq:waveMAP}), where
estimation is applied separately at each
resolution level $j$. The number of wavelet coefficients at the
$j$th resolution level is $n_j=2^j$. Let
$\pi_j(\kappa)>0,\;\kappa=0,\ldots,2^j$, be the prior distribution on the
number of non-zero wavelet coefficients, and let $\tau_j^2/n$ be
their level-wise prior variance, $j_0 \leq j < J$; the corresponding level-wise
variance ratios are $\gamma_j=\tau_j^2/\sigma^2$.

\begin{theorem} \label{th:level0}
Let $\psi$ be a mother wavelet of regularity $r$ and let
$\hat{f}_n(\cdot)$ be the corresponding level-wise wavelet
estimator (\ref{eq:waveMAP})  of $f$ in the Gaussian white noise model
(\ref{eq:FaGauWh}), where $f \in B^s_{p,q}(M)$, $0 < p,q \leq
\infty$, $1/p < s < r$ and $M>0$. Assume that there
exist positive constants $\gamma_{-}$ and $\gamma_{+}$
such that $\gamma_{-} \leq \gamma_j \leq \gamma_{+}$ for
all $j=j_0,\ldots,J-1$. Let the priors $\pi_j$ satisfy the
following conditions for all $j=j_0,\ldots,J-1$:
\begin{enumerate}
\item $\pi_j(0) \geq 2^{-c_1j}$ for some $c_1>0$;
\item $\pi_j(\kappa) \geq (\kappa2^{-j})^{c_2 \kappa}$ for all
$\kappa=1,2,\ldots,\alpha_j 2^j$, where $c_2>0$ and $0 < c_\alpha \leq
\alpha_j \leq \exp\{-c(\gamma_j)\}$ for some constant $c_\alpha > 0$,
and the function $c(\gamma_j)=8(\gamma_j+3/4)^2$ was defined in
Proposition \ref{prop:prop1};
\item $\pi_j(2^j) \geq e^{-c_0 2^j}$ for some $c_0>0$. 
\end{enumerate}
Then, as $n \rightarrow \infty$,
$$
\sup_{f \in B^s_{p,q}(M)} E(||\hat{f}_n-f||_2^2)=
O\Big(n^{-\frac{2s}{2s+1}}\Big).
$$
\end{theorem}

\medskip

For $f \in B^s_{p,q}(M)$, the sequence of its wavelet coefficients
at the $j$th resolution level belongs to $l_p[\eta_j]$, where
$\eta_j=C_0 n^{1/2} 2^{-j(s+1/2)}$ for some $C_0>0$
(Meyer, 1992, Section 6.10). The conditions on the prior in Theorem
\ref{th:level0} ensure that all the four statements of the
Proposition \ref{prop:prop1} simultaneously hold at all resolution
levels $j_0 \leq j < J$ with $\beta=0$, and one can exploit any of
them at each resolution level. It is necessary for adaptivity
of the resulting level-wise wavelet estimator (\ref{eq:waveMAP}).

%The assumptions of Theorem \ref{th:level0} are, in fact, not too
%restrictive. 
As we have mentioned in Section \ref{subsec:bounds}, all three
conditions of Theorem \ref{th:level0} hold, for example, for the truncated
geometric prior ${\rm TrGeom}(1-q_j)$, where $q_j$ are bounded away from zero 
and one.

It turns out that requiring a slightly more stringent condition on
$\pi_j(0)$, allows one also to estimate derivatives of $f$ by the
corresponding derivatives of its level-wise wavelet estimator
$\hat{f}_n$ at the optimal convergence rates. Such a plug-in
estimation of $f^{(m)}$ by $\hat{f}_n^{(m)}$ is, in fact, along the
lines of the vaguelette-wavelet decomposition approach of Abramovich
\& Silverman (1998).

Recall that, as $n \rightarrow \infty$, the minimax convergence rate for the 
$L^2$-risk of
estimating an $m$th derivative of the unknown response
function $f$ in the model (\ref{eq:FaGauWh}) over Besov balls
$B^{s}_{p,q}(M)$, where $0 \leq m < \min\{s,(s+1/2-1/p)p/2\}$, $0 < p,q \leq
\infty$ and $M>0$, is given by
$$
\inf_{\tilde{f}^{(m)}_n}\sup_{f \in B^s_{p,q}(M)}
E(||\tilde{f}^{(m)}_n-f^{(m)}||_2^2) \asymp n^{-2(s-m)/(2s+1)} 
$$
(Donoho et al., 1997; Johnstone and Silverman, 2005).

\medskip

The following Theorem \ref{th:level} is a generalization of Theorem
\ref{th:level0} for simultaneous level-wise wavelet estimation
of a function and its derivatives.

\begin{theorem} \label{th:level}
Let $\psi$ be a mother wavelet of regularity $r$ and let $\hat{f}_n$
be the level-wise wavelet estimator (\ref{eq:waveMAP}) of $f$ in the 
Gaussian white noise model (\ref{eq:FaGauWh}), where $f \in
B^s_{p,q}(M)$, $0 < p,q \leq \infty$, $1/p < s < r$
and $M>0$. Assume that there exist positive constants
$\gamma_{-}$ and $\gamma_{+}$ such that
$\gamma_{-} \leq \gamma_j \leq \gamma_{+}$ for all
$j=j_0,\ldots,J-1$. Let the priors $\pi_j$ satisfy the
following conditions for all $j=j_0,\ldots,J-1$:
\begin{enumerate}
\item $\pi_j(0) \geq 2^{-c_1j2^{-\beta j}}$ for some $\beta \geq 0$ and
$c_1>0$;
\item $\pi_j(\kappa) \geq (\kappa2^{-j})^{c_2 \kappa}$ for all
$\kappa=1,2,\ldots,\alpha_j 2^j$, where $c_2>0$ and $0 < c_\alpha \leq
\alpha_j \leq \exp\{-c(\gamma_j)\}$ for some constant $c_\alpha > 0$,
and the function $c(\gamma_j)=8(\gamma_j+3/4)^2$ was defined in
Proposition \ref{prop:prop1};
\item $\pi_j(2^j) \geq e^{-c_0 2^j}$ for some $c_0>0$. 
\end{enumerate}
Then, for all $m$th derivatives $f^{(m)}$ of $f$, where $ 0 \leq m
\leq \beta/2$ and $m < \min\{s,(s+1/2-1/p)p/2\}$, as $n \rightarrow \infty$,
$$
\sup_{f \in B^s_{p,q}(M)} E(||\hat{f}^{(m)}_n-f^{(m)}||_2^2)=
O\Big(n^{-\frac{2(s-m)}{2s+1}}\Big).   
$$
\end{theorem}

\medskip

Theorem \ref{th:level0} is evidently a particular case of Theorem
\ref{th:level} corresponding to the case $m=0$, for $\beta=0$ in the
condition on $\pi_j(0)$. Theorem \ref{th:level} shows that the same
proposed adaptive level-wise wavelet estimator (\ref{eq:waveMAP})
is simultaneously optimal
for estimating a function and an entire range of its derivatives.
This range is the same as that for the empirical
Bayes shrinkage and threshold estimators appearing in Theorem 1 of
Johnstone \& Silverman (2005). The proof of Theorem \ref{th:level}
is given in the Appendix.

\section{Numerical Study}\label{sec:sim}
\subsection{Preamble} \label{subsec:preamble}
In this section, we present a simulation study to illustrate the
performance of the developed level-wise wavelet estimator (\ref{eq:waveMAP})
and compare it with three empirical Bayes wavelet estimators:
the posterior mean and the posterior median 
of Johnstone \& Silverman (2005), and the Bayes Factor of
Pensky \& Sapatinas (2007); and two other estimators: the
block wavelet estimator NeighBlock of Cai \& Silverman (2001) and the 
complex-valued wavelet hard thresholding estimator  
of Barber \& Nason (2004). All the above Bayesian estimators and the block
wavelet estimator
are asymptotically minimax in a wide range of Besov balls.
Although no such theoretical results have been established so far for the 
complex-valued wavelet estimator,
it has performed well in simulations (Barber \& Nason, 2004).

In practice, one typically deals with 
discrete data of a sample size $n$ and the sampled data analog of the
Gaussian white noise model (\ref{eq:FaGauWh}) is the standard
nonparametric regression model
$$
Y_i=f(i/n)+\epsilon_i, \quad i=1,2,\ldots,n,
$$
where $\epsilon_i$ are independent $N(0,\sigma^2)$ random variables.
The corresponding global and level-wise Bayesian maximum a posteriori wavelet
estimation procedures then use the empirical wavelet coefficients
obtained by the discrete wavelet transforms of the data. However,
utilizing the machinery of Johnstone \& Silverman (2004a, 2005) for
development of appropriate boundary-corrected wavelet bases, one can
show that discretization does not affect the order of magnitude
of the accuracy of the resulting wavelet estimates (Johnstone \&
Silverman, 2004a, 2005; Petsa, 2009, Chapter 3).

The computational algorithms were
performed using the WaveLab and EbayesThresh software.
%http://www-stat.stanford.edu/software/software.html
%and the EbayesThresh software available at
%http://www-lmc.imag.fr/lmc-sms/Anestis.Antoniadis/EBayesThresh.
The entire study was carried out using the Matlab programming
environment.

\subsection{Estimation of parameters} \label{subsec:parameter}
To apply the proposed level-wise wavelet estimator (\ref{eq:waveMAP}) one 
should specify the priors $\pi_j$, the noise variance $\sigma^2$ and the prior
variances $\tau^2_j$ or, equivalently, the variance ratios
$\gamma_j=\tau_j^2/\sigma^2$. We used the truncated geometric priors
${\rm TrGeom}(1-q_j)$ discussed in Section \ref{subsec:level}. Since
the parameters $\sigma^2,\;q_j$ and $\gamma_j$ are rarely known 
a priori in practice, they should be estimated from the data in the
spirit of empirical Bayes.

The unknown $\sigma$ was robustly estimated by the median of the
absolute deviation of the empirical wavelet coefficients at the
finest resolution level $J-1$, divided by 0.6745 as suggested by
Donoho \& Johnstone (1994a), and usually applied in practice. For a
given $\sigma$, we then estimate $q_j$ and $\gamma_j$ by the
conditional likelihood approach of Clyde \& George (1999).

Consider the prior model described in Section \ref{subsec:test}. The
corresponding marginal likelihood of the observed empirical wavelet
coefficients, say $Y_{jk}$, at the $j$th resolution level is then
given by
$$
L(q_j,\gamma_j;Y_j) \propto \sum_{\kappa=0}^{2^j}\pi_j(\kappa) {2^j
\choose \kappa}^{-1} (1+\gamma_j)^{-\kappa/2} \sum_{x_i: \sum_k
x_{ik}=\kappa} \exp\Bigg\{\frac{\gamma_j \sum_k
x_{ik}Y_{jk}^2}{2\sigma^2(1+\gamma_j)}\Bigg\},
$$
where $\pi_j(\kappa)=(1-q_j)q_j^{\kappa}/(1-q_j^{2^j+1})$ and $x_i$
are indicator vectors. Instead of direct maximization of
$L(q_j,\gamma_j;Y_j)$ with respect to $q_j$ and $\gamma_j$, regard
the indicator vector $x$ as a latent variable and consider the
corresponding log-likelihood for the augmented data $(Y_j,x)$,
i.e.,
\begin{equation}
l(q_j,\gamma_j;Y_j,x)= c+\log\pi_j(\kappa)-\log{2^j
\choose \kappa} -\frac{\kappa}{2}\log(1+\gamma_j) + \frac{\gamma_j \sum_k
x_{ik}Y_{jk}^2}{2\sigma^2(1+\gamma_j)}, \label{eq:loglike}
\end{equation}
where $c$ is a constant.
The EM-algorithm iteratively alternates between computation of the
expectation of $l(q_j,\gamma_j;Y_j,x)$ in (\ref{eq:loglike}) with
respect to the distribution of $x$ given $Y_j$ evaluated using the
current estimates for the parameters' values at the E-step, and
updating then the parameters by maximizing it with respect to $q_j$ and
$\gamma_j$ at the M-step. However, for a general prior distribution
$\pi_n$ and for the truncated geometric prior, in particular, the
EM-algorithm does not allow one to achieve analytic expressions on the
E-step. Instead, we apply the conditional likelihood estimation
approach originated by George \& Foster (2000) and adapted to the
wavelet estimation context by Clyde \& George (1999). The approach
is based on evaluating the augmented log-likelihood
(\ref{eq:loglike}) at the mode for the indicator vector $x$ at
the E-step rather than using the mean as in the original
EM-algorithm (Abramovich \& Angelini, 2006).

For a fixed number $\kappa$ of its non-zero entries, it is evident
from (\ref{eq:loglike}) that the most likely vector
$\hat{x}(\kappa)$ is $\hat{x}_i(\kappa)=1$ for the $\kappa$ largest
$|Y_{jk}|$ and zero otherwise. For the given $\kappa$, maximizing
(\ref{eq:loglike}) with respect to $\gamma_j$ after some algebra yields
$\hat{\gamma}_j(\kappa)=\max\Big\{0,\sum_{k=1}^\kappa Y^2_{(k)}/
(\kappa \sigma^2)-1\Big\}$. To simplify maximization with respect to
$q_j$, approximate the truncated geometric distribution $\pi_j$ in
(\ref{eq:loglike}) by a non-truncated one. This approximation does
not strongly affect the results, especially at sufficiently high
resolution levels, and allows one to obtain analytic solutions for
$\hat{q}_j$, i.e., $\hat{q}_j(\kappa)=\kappa/(\kappa+1)$. It is now
straightforward to find $\hat{\kappa}$ that maximizes
(\ref{eq:loglike}) together with the corresponding
$\hat{\gamma}_j(\hat{\kappa})$ and $\hat{q}_j(\hat{\kappa})$. The
above conditional likelihood approach results therefore in rapidly
computable estimates for $\gamma_j$ and $q_j$ in closed forms.

\subsection{Simulation study}
\label{subsec:simFF}
We now present and discuss the results of the simulation study.
For all three empirical Bayes wavelet estimators, we used the
Double-exponential prior, where the corresponding
prior parameters were estimated level-by-level by
marginal likelihood maximization, as described in Johnstone \&
Silverman (2005). The prior parameters for the proposed level-wise 
wavelet estimator (\ref{eq:waveMAP}) were estimated by conditional likelihood
maximization procedure described in Section \ref{subsec:parameter} above. For
the block wavelet estimator, the lengths of the blocks and the thresholds
were selected as suggested
by Cai \& Silverman (2001). Finally, for all competing methods,
$\sigma$ was estimated by the median of the absolute value of the
empirical wavelet coefficients at the finest resolution level
divided by 0.6745 as discussed in Section \ref{subsec:parameter}.

In the simulation study, we evaluated the above six wavelet estimators
for a series of test functions.
We present here the results for the nowadays standard Bumps, Blocks, Doppler and
Heavisine functions of Donoho \& Johnstone (1994a), 
and Wave (Marron et al., 1998; Antoniadis et al., 2001) and Peak (Angelini et al., 2003) functions defined, respectively, as
\begin{equation*}f(t)=0.5 + 0.2\cos{4\pi t}+ 0.1\cos{24\pi
t}, \;\;\; 0 \leq t \leq 1 \end{equation*}
%\begin{equation*}
%f(t)= \left(0.32 + 0.6t+0.3 \exp{-100(t-0.3)^2}\right)\II_{[0,0.8]}(t)
%+ \left( -0.28 + 0.6 t + 0.3
%\exp{-100(t-1.3)^2}\right)\II_{(0.8,1]}(t)
%\end{equation*}
and
\begin{equation*}f(t)=\exp\{-|t-0.5|\}, \;\;\;0 \leq t \leq 1 .\end{equation*}

\begin{figure}
\begin{center}
\includegraphics[width=1.0\textwidth]{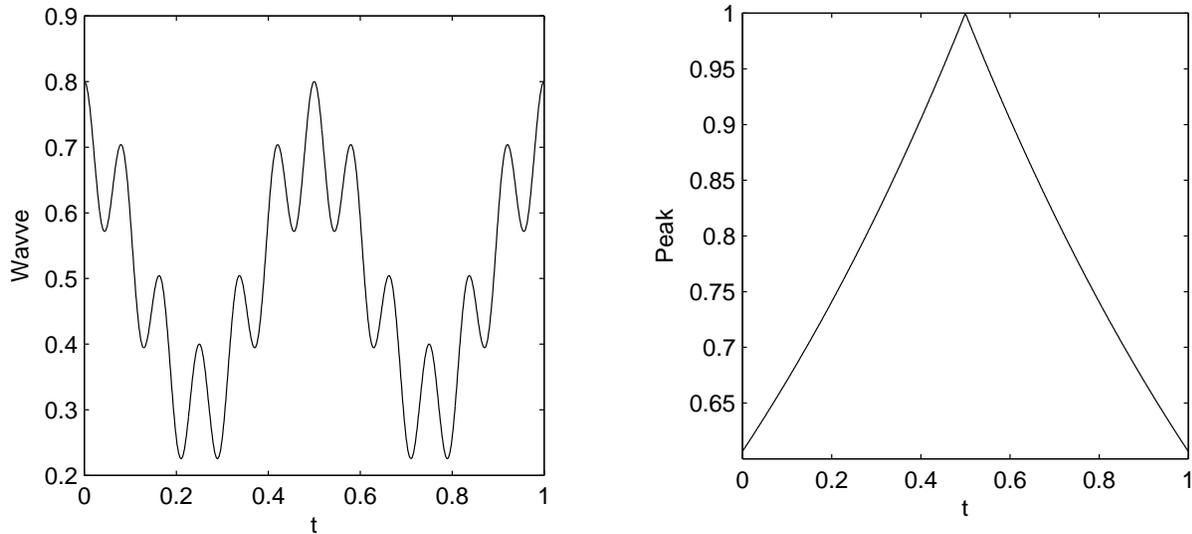}
\caption{\label{testfunctions}Wave (left) and Peak (right) test functions}
\end{center}
\end{figure}
\noindent
See Figure \ref{testfunctions} for Wave and Peak test functions.

For each test function, $M=100$ samples were generated by adding independent
Gaussian noise $\varepsilon \sim N(0,\sigma^2)$ to $n=$ 256, 512 and
1024 equally spaced points on [0,1].
The value of the root signal-to-noise ratio
was taken to be 3, 5 and 7 corresponding respectively to high, moderate and low
noise levels.
%$$
%{\rm SNR}(f,\sigma) =  \sigma^{-1}\
%\left(\frac{1}{n}\sum_{i=1}^{n}(f(t_i)-\bar{f})^2\right)^{1/2},\;\;\;
%\bar{f} = \frac{1}{n}\sum_{i=1}^{n}f(t_i).
%$$
The goodness-of-fit for an estimator $\hat{f}$ of $f$ in a single replication was measured
by its mean squared error. 
%defined as
%$$
%{\rm MSE}(f,\hat{f}) = \frac{1}{n} \sum_{i=1}^{n}(\hat{f}(t_i)-f(t_i))^2.
%$$

For brevity, we report the results only for  $n=1024$ using the
compactly supported mother wavelet Coiflet 3 (Daubechies, 1992, p.258)
and the Lawton mother wavelet (Lawton, 1993) for the complex-valued wavelet 
estimator. The primary resolution level was $j_0=4$. Different choices of sample
sizes and wavelet functions basically yielded similar results in magnitude. 

The sample distributions of
mean squared errors over replications for different wavelet estimators in the 
conducted simulation study were
typically asymmetrical and affected by outliers. Therefore,
we preferred the sampled medians of mean squared errors rather than means to gauge
the estimators' goodness-of-fit. Thus, for each wavelet estimator, test function
and noise level, we calculated the sample median of mean squared errors over all 
100 replications. To quantify the comparison between the competing wavelet 
estimators over various test functions and noise levels, for each model
we found the best wavelet estimator among the six, i.e., the one
achieving the minimum median mean squared error. We then evaluated the relative
median mean squared error of each estimator defined as the 
ratio between the minimum and the estimator's median mean squared errors;
%$\min_{1 \leq j \leq 6}\{\text{Median}(\text{MSE}_j)\}/\text{Median}(\text{MSE}_i)$,
%$i=1,\ldots,6$ 
see Table \ref{median}. 

As expected, Table
\ref{median} shows that there is no uniformly best
wavelet estimator. Each one has its own favorite and challenging cases,
and its relative performance strongly depends on a specific test function.
Thus, the complex-valued estimator indeed demonstrates excellent results for 
Donoho \& Johnstone's functions as it has been reported in Barber \& Nason (2004),
but is much less successful for Peak and Wave.
The block estimator is the best for the Peak and Doppler but the
worst for Blocks and Bumps. The proposed Bayesian estimator (\ref{eq:waveMAP})
outperforms others for 
Wave but is less efficient for Donoho \& Johnstone's (1994a) examples.
Interestingly, the relative performance of the estimators is much less sensitive
to the noise level. For each of the test functions, the corresponding best estimator is 
essentially the same for all noise levels. 

The minimal relative median of 
mean squared errors of an estimator over all cases 
reflects its inefficiency
at the most challenging combination of a test function and noise level, and can
be viewed as a natural measure of its robustness. In this sense,
the posterior mean estimator is the most robust although it is not
the winner in any particular case. 
%However, even its minimal relative median is below 0.65.

\begin{table}
\caption{Relative median mean squared errors (MSE) for various
test functions, 
levels of the root signal-to-noise ratio (RSNR) and different wavelet 
estimators: the proposed Bayesian estimator (MAP), Bayes Factor (BF),
posterior median (Postmed), posterior mean (Postmean), block (Block)
and complex-valued hard thresholding (CW).}{%
%The minimal relative median MSE for each estimator is bold.}
\begin{tabular}{lccccccc}
 signal & RSNR & MAP & BF & Postmed & Postmean & Block & CW \\
Peak &3 &0.8697& 0.1763 & 0.8279&0.6589&1&0.5795  \\
     &5&  0.7772 & 0.1497&0.7864&0.6525&1&0.6234 \\
     &7&  0.8033& 0.186&0.8501&0.6958&1&0.6979 \\
 \\
Wave &3&   1&0.5614 &0.9841& 0.9103&0.4570&0.9189 \\
     &5&  0.9841&0.4603 & 1 &0.9165&0.6072 &0.8265\\
     &7&  1& 0.6241&0.9900&0.9303&0.7498&0.7793\\
 \\
Bumps &3&  0.5968 &0.6254& 0.6814&0.7569&0.4769&1 \\
      &5&  0.5221 & 0.5641 &0.5893&0.6671&0.4788&1 \\
      &7&  0.5132 & 0.5537 &0.5707& 0.6420& 0.5202&1\\
 \\
Blocks &3&  0.6595 &0.6807& 0.8815&0.9500&0.5606&1 \\
       &5&  0.6875&0.727 & 0.8541 &0.9065&0.4416&1 \\
       &7&  0.6921& 0.7134&0.7806&0.8535&0.4288&1\\
 \\
Doppler &3&  0.7214 &0.611& 0.8277&0.8709&0.9878 &1\\
       &5&  0.6962&0.6739 & 0.8116 &0.8583&1&0.9119 \\
       &7&  0.7655&0.7122&0.8236&0.883&1&0.9382\\
 \\
HeaviSine &3&  0.7523 &0.3566& 0.9333&0.9154&0.8406&1 \\
       &5&  0.6640&0.3764 & 0.8622 &0.8427&0.5796 &1\\
       &7&  0.6931& 0.3505&0.8298&0.8424&0.5028&1\\
\end{tabular}} \label{median}
\end{table}

As a ``by-product'', we also compared different thresholding estimators in terms of
sparsity measured by the average percentage of non-zero wavelet coefficients
which remained after thresholding; see Table
\ref{sparse}. The posterior mean estimator was not
included in this comparison since it is a non-linear shrinkage but not 
a thresholding estimator. Similar to the previous results on the goodness-of-fit,
the relative sparsity strongly depends on the test function.
However, except for the Doppler
example, the estimator (\ref{eq:waveMAP}) is consistently the most 
sparse among Bayesian estimators.

%Finally, it is important to note that
%the MAP testimator was faster
%(in terms of CPU time) than the PostMean, PostMed and BF estimators
%but, as expected, slower than the Block.

\begin{table}
\caption{Average percentages of remaining coefficients 
for various test functions, levels of the root signal-to-noise ratio 
and different wavelet thresholding estimators.}{%
\begin{tabular}{lcccccc}
 signal & RSNR & MAP & BF & Postmed & Block & CW \\
Peak &3& 7.57 &14.87& 8.92& 1.64& 1.74 \\
     &5& 5.62 &15.89& 8.39& 1.61& 1.93 \\
     &7& 4.26 &12.93& 7.81& 1.63& 2.05 \\
     \\
Wave &3& 11.39 &18.81&12.79 & 5.21& 5.52 \\
     &5& 11.51&20.19&12.52 & 6.30&  6.28 \\
     &7& 10.38& 20.84& 13.07& 6.30& 7.01 \\
     \\
Bumps &3& 10.63 & 12.13 &10.86& 16.63& 12.23 \\
      &5& 11.17 & 12.60 &12.45& 21.13& 14.52 \\
      &7& 12.65 & 13.90 &13.80& 23.70& 16.03 \\
      \\
Blocks &3& 17.10& 15.32& 11.62& 12.27& 8.39 \\
       &5& 10.39& 12.20& 11.72& 18.03&12.07 \\
       &7& 11.12& 12.87& 12.63& 22.47&14.20 \\
       \\
Doppler &3&11.46& 14.73&  8.58&  5.69& 5.13 \\
        &5& 7.24&  9.23&  6.52&  6.63& 6.60 \\
        &7& 6.42&  9.58&  6.57&  7.27& 7.86 \\
	\\
HeaviSine &3& 6.35&19.27&10.75& 1.98& 2.17 \\
          &5& 8.11&19.73&10.62& 3.17& 2.69 \\
          &7&10.87&18.52&12.55& 4.00& 3.39 \\
\end{tabular}} \label{sparse}
\end{table}

Apart from providing a theoretical justification, the presented numerical
results show that the proposed estimator demonstrates good performance in finite 
sample settings and can, therefore, be viewed as a contribution to the 
list of useful wavelet-based function estimation tools.

%\section{Concluding Remarks } \label{sec:conclusions}
%Although we considered only quadratic losses in our exposition, we
%believe that the obtained results can be extended to more general
%losses similar to those in Donoho \& Johnstone (1994b) and Johnstone
%\& Silverman (2004b, 2005). Furthermore, the proposed methodology
%can be adapted to derive pointwise optimal (in the minimax sense)
%MAP wavelet testimators of the unknown response function and its
%derivatives in the standard Gaussian white noise model, as in Cai
%(2003). However, appropriate adjustments are needed for each
%specific problem at hand, and we hope to address these issues
%elsewhere.

\section*{Acknowledgments}
Felix Abramovich and Vadim Grinshtein were supported by the Israel
Science Foundation grant ISF-248/08. The authors would like to
thank the Editor, the Associate Editor and the three anonymous
referees for helpful comments on improvements to this paper.
Valuable remarks of Bernice Oberman are gratefully acknowledged.

\section{Appendix} \label{sec:appendix}
Throughout the proofs we use $C$ to denote a generic positive constant, not
necessarily the same each time it is used, even within a single equation.

\subsection{Proof of Proposition \ref{prop:prop1}} \label{subsec:prop}
We start the proof of the proposition with the following lemma that
establishes the bounds for binomial coefficients:

\begin{lemma} \label{lem:lem1}
For all $n \geq 2$ and $\kappa=1,2,\ldots,n-1$,
\begin{equation}
\left(\frac{n}{\kappa}\right)^\kappa \leq {n \choose \kappa} <
\left(\frac{ne}{\kappa}\right)^{\kappa}. \label{eq:lem11}
\end{equation}
In particular, for $\kappa \leq n/e$,
\begin{equation}
{n \choose \kappa} < \left(\frac{n}{\kappa}\right)^{2 \kappa}.
\label{eq:lem12}
\end{equation}
\end{lemma}

This lemma generalizes Lemma A.1 of Abramovich et al. (2007), 
where the upper bound similar to that in
(\ref{eq:lem11}) was obtained for $\kappa=o(n)$.

{\sc Proof of Lemma \ref{lem:lem1}}. The obvious lower bound for the
binomial coefficient in (\ref{eq:lem11}) has been shown in Lemma A.1 of Abramovich et al. (2007). 
To prove the upper bound in
(\ref{eq:lem11}), note that using Stirling's formula one has
\begin{equation}
{n \choose \kappa} \leq \left(\frac{n}{e}\right)^n
\left(\frac{e}{n-\kappa}\right)^{n-\kappa}
\left(\frac{e}{\kappa}\right)^\kappa =
\left(\frac{n}{\kappa}\right)^\kappa
\left(\frac{n}{n-\kappa}\right)^{n-\kappa} \label{eq:stir}
\end{equation}
for all $n \geq 2$ and $\kappa=1,2,\ldots,n-1$.

Note that $\log\{x/(x-1)\} < 1/(x-1)$ for all $x>1$. In particular, for
$x=n/\kappa$ it implies $\log\{n/(n-\kappa)\} < \kappa/(n-\kappa)$ and,
therefore, $(\frac{n}{n-\kappa})^{n-\kappa} < \exp(\kappa)$ that
together with (\ref{eq:stir}) completes the proof of
(\ref{eq:lem11}).

The second statement (\ref{eq:lem12}) of the lemma is an immediate
consequence of (\ref{eq:lem11}) for $\kappa \leq n/e$. This
completes the proof of Lemma~\ref{lem:lem1}.

\medskip

We now return to the proof of Proposition \ref{prop:prop1} and
consider separately all the four cases covered by the proposition.
The proof will exploit the general results of Abramovich et al. (2007) 
on the upper error bounds for the $l^2$-risk of the
estimator (\ref{eq:MAP}) adapting them also for the case where the variance in the
Gaussian sequence model $(\ref{eq:sequence})$ may depend on $n$.

\medskip
\noindent Case 1. Under the condition $\pi_n \geq e^{-c_0 n},\;c_0>0$, 
the definition (\ref{eq:MAP}) of $\hat{\mu}$ and (\ref{eq:penalty}) 
immediately imply
$
||y-\hat{\mu}||^2 \leq ||y-\hat{\mu}||^2 + P_n(\hat{k}) \leq  P_n(n) = O(n\sigma_n^2).
$
Thus,
$$
E(||\hat{\mu}-\mu||^2) \leq 2\{E(||y-\hat{\mu}||^2) + E(||y-\mu||^2)\} = O(n\sigma_n^2).
$$

\medskip
\noindent Case 2. Applying Corollary 1 of Abramovich et al. (2007)
for $\kappa=0$ yields
$$
E(||\hat{\mu}-\mu||^2_2) \leq c_0(\gamma_n)\Big\{\sum_{i=1}^n
\mu_i^2+2\sigma_n^2(1+1/\gamma_n)\log\pi_n^{-1}(0)\Big\}
+c_1(\gamma_n)\{1-\pi_n(0)\}\sigma_n^2,
$$
where the exact expressions for $c_0(\gamma_n)$ and $c_1(\gamma_n)$ are given
in Theorem 2 of Birg\'e \& Massart (2001) with their $K=1+1/(2\gamma_n)$;
see the proof of Theorem 1 of Abramovich et al. (2007). In particular,
under the assumptions of the proposition on the boundness of $\gamma_n$,
the functions $c_0(\gamma_n)$ and $c_1(\gamma_n)$ are also bounded from above.
For $2 \leq p \leq \infty$, the
least favorable sequence $\mu_0$ that maximizes $\sum_{i=1}^n
\mu_i^2$ over $l_p[\eta_n]$ is $\mu_{01}=\ldots=\mu_{0n}=C_n
n^{-1/p}=\eta_n \sigma_n$. As $n \rightarrow \infty$, one then has
\begin{align*}
E(||\hat{\mu}-\mu||^2_2) & \leq c_0(\gamma_n)\left\{\sigma_n^2
\eta^2_n n + 2 \sigma_n^2 (1+1/\gamma_n) \log \pi_n^{-1}(0)\right\}+
c_1(\gamma_n)\{1-\pi_n(0)\}\sigma_n^2\\
& = O(\sigma_n^2 n \eta_n^2)+O(\sigma_n^2 n^{-\beta}\log n).
\end{align*}

\medskip
\noindent Case 3. This is essentially a sparse case
considered in
Abramovich et al. (2007) and its proof is a direct consequence of their Theorem 6.

\medskip
\noindent Case 4. The proof for this case is similar to that
of Case 2 except that for $0<p<2$, the least favorable sequences
$\mu_0$ that maximize $\sum_{i=1}^n \mu_i^2$ over $\mu \in
l_p[\eta_n]$ are permutations of the spike $(C_n,0,\ldots,0)$
and therefore $\sum_{i=1}^n \mu_{0i}^2 \leq \sigma_n^2 n^{2/p}
\eta^2_n $. Repeating the arguments used in the proof of Case 2 for
$\kappa=0$, under the requirements of the proposition on boundedness
of $\gamma_n$, we then get as $n \rightarrow \infty$,
\begin{align*}
E(||\hat{\mu}-\mu||^2_2) &\leq c_0(\gamma_n)\left\{\sigma_n^2 n^{2/p}
\eta^2_n + 2(1+1/\gamma_n)\sigma_n^2 \log \pi^{-1}_n(0)\right\}+
c_1(\gamma_n)\{1-\pi_n(0)\}\sigma^2_n
\\&= O(\sigma_n^2 n^{2/p} \eta^2_n)+O(\sigma_n^2 n^{-\beta}\log n).
\end{align*}
for all $\eta^p_n < n^{-1}(2 \log n)^{p/2}$. This completes the proof
of Theorem~\ref{prop:prop1}.

\subsection{Proof of Theorem \ref{th:global}} \label{subsec:proofglobal}
Let $R_j=\sum_{k=0}^{2^j-1}E\{(\hat{\theta}_{jk}-\theta_{jk})^2\},\;j
\geq j_0-1$, be the $L^2$-risk of the global wavelet estimator (\ref{eq:waveMAP})
at the $j$th resolution level. Due to the Parseval relation,
$E(||\hat{f}_n-f||^2)=\sum_{j \geq j_0-1} R_j$. Scaling coefficients
are not thresholded and therefore $R_{j_0-1}=C_{j_0}\sigma^2
n^{-1}=o(n^{-2s/(2s+1)})$ as $n \rightarrow \infty$. At very high
resolution levels, where $j \geq J$, all wavelet coefficients
$\hat{\theta}_{jk}$ are set to zero and, therefore, as $n \rightarrow \infty$,
$$
\sum_{j=J}^\infty R_j=\sum_{j=J}^\infty \sum_{k=0}^{2^j-1}
\theta_{jk}^2= O(n^{-2s'})=o(n^{-2s/(2s+1)}),
$$
where $s'=s+1/2-1/\min(p,2)$ (Johnstone \& Silverman, 2005).

Consider now $\sum_{j=j_0}^{J-1}R_j$. The set of wavelet
coefficients $\{\theta_{i}\}$ of a function $f \in B^s_{p,q}(M)$
lies within a weak $l_{r}$-ball of a radius $a M$ with $r=2/(2s+1)$,
where the constant $a$ depends only on a chosen wavelet basis:
$m_r[\eta_n]=\{\theta: |\theta|_{(i)} \leq (a M)i^{-1/r}\}$ 
(Donoho, 1993, Lemma 2). The corresponding normalized radius
$\eta_n=(\sigma/\surd{n})^{-1}\tilde{n}^{-1/r}a M=O(n^{-s})$, where
$\tilde{n}=n-2^{j_0} \sim n$ for large $n$.

Under the conditions of the theorem, one can then apply Theorem 6 of
Abramovich et al. (2007) for $m_r[\eta_n]$ to get 
$$
\sum_{j=j_0}^{J-1}R_j \leq \sup_{\theta \in
m_{r}[\eta_n]}E(||\hat{\theta}-\theta||^2_2) =O\left\{\eta_n^{r}(2\log
\eta_n^{-r})^{1-r/2}\right\} = O\left\{\left(\frac{\log
n}{n}\right)^{2s/(2s+1)}\right\} 
$$
as $n \rightarrow \infty$.
This completes the proof of Theorem~\ref{th:global}.

\subsection{Proof of Theorem \ref{th:level}} \label{subsec:prooflevel}
Let $R_j=\sum_{k=0}^{2^j-1}E(\hat{\theta}_{jk}-\theta_{jk})^2,\;j
\geq j_0-1$, be now the $L^2$-risk of the level-wise version of the wavelet
estimator (\ref{eq:waveMAP}) at the $j$th resolution level. 
Johnstone \& Silverman
(2005, Section 5.6) showed that $E(||\hat{f}_n^{(m)}-f^{(m)}||^2)
\asymp \sum_{j \geq j_0-1} 2^{2mj}R_j$.

For any $f \in B^s_{p,q}(M)$, the sequence of its wavelet
coefficients at the $j$th resolution level belongs to a strong $l_p$-ball
of a normalized radius $\eta_j=C_0 n^{1/2}2^{-j(s+1/2)}$ for some $C_0>0$ 
(Meyer, 1992, Section 6.10).

Define
$$
j_1=\frac{1}{2s+1}\log_2 \left(\frac{n C_0^2}{c_\alpha^{2/p}}\right)
\sim \frac{1}{2s+1} \log_2 n.
$$
For sufficiently large $n$, $j_1>j_0$. Note that $\eta^p_j \geq
c_\alpha$ for $j \leq j_1$ and $\eta^p_j < c_\alpha$ for $j > j_1$
with obvious modifications for $p=\infty$. Consider the following
cases:

\noindent
\newline 1. Scaling coefficients: $j=j_0-1$. Similarly to the
global wavelet estimator, for a fixed primary resolution level
$j_0$, $2^{2m(j_0-1)}R_{j_0-1}=O(n^{-1})=o(n^{-2(s-m)/(2s+1)})$ as
$n \rightarrow \infty$.

\noindent
\newline 2. Coarse resolution levels: $j_0 \leq j \leq j_1$. Applying the first
statement of Proposition \ref{prop:prop1} for each level one has
$$
\sum_{j=j_0}^{j_1}2^{2mj}R_j \leq C \sum_{j=j_0}^{j_1}
2^{2mj}n^{-1}\sigma^2 n_j  \leq C n^{-1} \sum_{j=j_0}^{j_1}
2^{(2m+1)j} = O\Big(n^{-2(s-m)/(2s+1)}\Big)
$$
as $n \rightarrow \infty$.

\noindent
\newline 3. Middle and high resolution levels: $j_1 < j < J$.
Consider separately the cases (a) $2 \leq p \leq \infty$ and
(b) $0 < p < 2$.

\noindent
\newline (a) $2 \leq p \leq \infty$. Under the conditions of the theorem,
the second statement of Proposition \ref{prop:prop1} at the $j$th
resolution level yields
$$
R_j \leq C n^{-1} \big(n_j \eta_j^2 + n_j^{-\beta} \log n_j\big) \leq
C\big(2^{-2js}+n^{-1}2^{-\beta j} j \big)
$$
and, hence, as $n \rightarrow \infty$,
\begin{align*}
\sum_{j=j_1+1}^{J-1}2^{2mj}R_j \leq C
\Big(2^{-2j_1(s-m)}+n^{-1}J^2\Big) &\leq
C\Big(n^{-2(s-m)/(2s+1)}+n^{-1}\log_2^2 n\Big) \\ &=
O\Big(n^{-2(s-m)/(2s+1)}\Big).
\end{align*}

\noindent
\newline (b) $0 < p < 2$.
%Define
%$$
%j_2=\frac{1}{2(s+1/2-1/p)} \log_2\left(\frac{n}{\log_2 n}\right)
%$$
%Obviously, $j_1 < j_2 < J$ for sufficiently large $n$.
Let $j_2$ be the largest integer for which $\eta_j^p \geq
n_j^{-1}(2\log n_j)^{p/2}$. One can easily verify that $j_1 < j_2 <
J$.

Using the monotonicity arguments, $\eta^p_j \geq n_j^{-1} (2\log
n_j)^{p/2}$ for all middle resolution levels $j_1 < j \leq j_2$.
One can then apply the third statement of Proposition
\ref{prop:prop1}, and after some algebra, to get, for
$m<(s+1/2-1/p)p/2$,
\begin{align*}
\sum_{j=j_1+1}^{j_2}2^{2mj}R_j & \leq  C n^{-1} \sum_{j=j_1+1}^{j_2}
2^{(2m+1)j} n^{p/2} 2^{-jp(s+1/2)} \left\{\log(n^{-p/2}2^{jp(s+1/2)})\right\}^{1-p/2} \\
& \leq   C
n^{-(1-p/2)}2^{-j_1p(s+1/2-(2m+1)/p)}\log\left(n^{-p/2}2^{j_1p(s+1/2)}\right) \\
& =O\Big(n^{-2(s-m)/(2s+1)}\Big)  
\end{align*}
as $n \rightarrow \infty$.

At high resolution levels $j_2 < j < J$, $\eta^p_j < n_j^{-1} (2\log
n_j)^{p/2}$, and the fourth statement of Proposition
\ref{prop:prop1} implies
$$
R_j \leq C\big(2^{-2j(s+1/2-1/p)}+n^{-1}2^{-j\beta}j\big).
$$
Hence, for $0 \leq m \leq \beta/2$ and $m<\min\{s,(s+1/2-1/p)p/2\}$,
one has
$$
\sum_{j=j_2+1}^{J-1}2^{2mj}R_j \leq C
\big(2^{-2(j_2+1)(s+1/2-1/p-m)}+n^{-1}J^2\big) = S_1 + S_2,
$$
where evidently
$S_2=O(n^{-1}\log^2_2n)=o\big(n^{-2(s-m)/(2s+1)}\big)$ as $n
\rightarrow \infty$. From the definition of $j_2$,
$2^{(j_2+1)(s+1/2-1/p)} > \surd\{nC/(j_2+1)\} > \surd{(nC/\log_2n)}$,
which after some algebra yields $S_1=o\big(n^{-2(s-m)/(2s+1)}\big)$
as $n \rightarrow \infty$.

\noindent
\newline 4. Very high resolution levels: $j \geq J$. Using the results of
Johnstone \& Silverman (2005), as $n \rightarrow \infty$, the tailed sum
$$
\sum_{j \geq J}
2^{2mj}R_j=O\big(n^{-2(s'-m)}\big)=o\big(n^{-2(s-m)/(2s+1)}\big),
$$
where $s'=s+1/2-1/\min(p,2)$. Summarizing,
$$
\sum_{j \geq j_0-1}
2^{2mj}R_j = O\big(n^{-2(s-m)/(2s+1)}\big) 
$$
as $n \rightarrow \infty$.
This completes the proof of Theorem~\ref{th:level}.

\end{document}